\documentclass[12pt]{amsart}
\usepackage{amsmath, amssymb, latexsym, amsthm}
\usepackage[all]{xy}

\newcommand{\Bc}[9]{\bibitem{#1} {#2}, \emph{#3}, in: \textbf{#4} (#5), #6 #7, #8--#9.}
\newcommand{\arx}[1]{\texttt{http://arxiv.org/abs/#1}}
\newcommand{\CH}{the Continuum Hypothesis}
\newcommand{\BOfat}{\B_{\Omega_{\mathrm{fat}}}}%\B_{\mathrm{fat}}}
\newcommand{\inv}{^{-1}}
\newcommand{\Cantor}{{\{0,1\}^\N}}
\newcommand{\bP}{\mathbf{P}}
\newcommand{\fU}{\mathfrak{U}}
\newcommand{\fV}{\mathfrak{V}}
\newcommand{\fW}{\mathfrak{W}}
\newcommand{\sr}[2]{{\txt{$#1$\\$#2$}}}
\newcommand{\seq}[1]{\{#1\}_{n\in\N}}

\newcommand{\cI}{\mathcal{I}}
\newcommand{\A}{\forall}
\newcommand{\B}{\mathcal{B}}

\newcommand{\BG}{\B_\Gamma}
\newcommand{\BL}{\B_\Lambda}
\newcommand{\BT}{\B_\Tau}
\newcommand{\BTstar}{\B_{\Tau^*}}
\newcommand{\BO}{\B_\Omega}
\newcommand{\CG}{C_\Gamma}
\newcommand{\CL}{C_\Lambda}
\newcommand{\CT}{C_\Tau}
\newcommand{\CTstar}{C_{\Tau^*}}
\newcommand{\CO}{C_\Omega}

\newcommand{\Tau}{\mathrm{T}}

\newcommand{\cF}{\mathcal{F}}

\newcommand{\M}{\mathcal{M}}
\newcommand{\N}{\naturals}
\newcommand{\NN}{{{}^{\naturals}\naturals}}
\renewcommand{\inf}{P_\oo(\N)}

\newcommand{\R}{\reals}

\newcommand{\cU}{\mathcal{U}}
\newcommand{\Union}{\bigcup}
\newcommand{\cV}{\mathcal{V}}
\newcommand{\cW}{\mathcal{W}}

\newcommand{\Impl}{\Rightarrow}
\long\def\forget#1\forgotten{}

\renewcommand{\b}{\mathfrak{b}}
\renewcommand{\t}{\mathfrak{t}}

\renewcommand{\c}{\mathfrak{c}}

\renewcommand{\i}{\item}

\newcommand{\oo}{\infty}
\renewcommand{\r}{\mathfrak{r}}
\renewcommand{\u}{\mathfrak{u}}
\newcommand{\h}{\mathfrak{h}}
\newcommand{\p}{\mathfrak{p}}

\newcommand{\w}{\omega}
\newcommand{\x}{\times}

\newcommand{\nin}{\not\in}

\newcommand{\sbst}{\subseteq}
\newcommand{\spst}{\supseteq}
\newcommand{\sm}{\setminus}

\newcommand{\as}{\subseteq^*}%{\let\proclaim\relax}
\renewcommand{\pi}{pseudo-intersection}
\renewcommand{\|}{\upharpoonright}%\restriction}

\renewcommand{\>}{\rangle}
\newcommand{\E}{\exists}

\newcommand{\cov}{\mathsf{cov}}
\newcommand{\add}{\mathsf{add}}

\newcommand{\non}{\mathsf{non}}

\newcommand{\impl}{\to}

\newtheorem{thm}{Theorem}[section]
\newtheorem{prop}[thm]{Proposition}
\newtheorem{prob}[thm]{Problem}
\newtheorem{lem}[thm]{Lemma}
\newtheorem{cor}[thm]{Corollary}
\newtheorem{conj}[thm]{Conjecture}
\theoremstyle{definition}

\theoremstyle{remark}

\newcommand{\be}{\begin{enumerate}}
\newcommand{\ee}{\end{enumerate}}
\newcommand{\bi}{\begin{itemize}}
\newcommand{\ei}{\end{itemize}}
\newcommand{\sone}{\mathsf{S}_1}

\renewcommand{\split}{\mathsf{Split}}

\newcommand{\naturals}{{\mathbb N}}
\newcommand{\reals}{{\mathbb R}}

\author{Boaz Tsaban}
\thanks{Partially supported by the Golda Meir Fund and
the Edmund Landau Center for Research in Mathematical Analysis and Related Areas,
sponsored by the Minerva Foundation (Germany).}
\address{Einstein Institute of Mathematics, Hebrew University of Jerusalem,
Givat Ram, Jerusalem 91904, Israel}
\email{tsaban@math.huji.ac.il}
\urladdr{http://www.cs.biu.ac.il/\~{}tsaban}
\title{The combinatorics of splittability}

\begin{document}
\begin{abstract}
Marion Scheepers, in his studies of the combinatorics of open covers,
introduced the property $\split(\fU,\fV)$ asserting that a cover of
type $\fU$ can be split into two covers of type $\fV$.
In the first part of this paper
we give an almost complete classification of all properties of this
form where $\fU$ and $\fV$
are significant families of covers which appear in the literature
(namely, large covers, $\omega$-covers, $\tau$-covers, and $\gamma$-covers),
using combinatorial characterizations of these properties in terms related
to ultrafilters on $\N$.

In the second part of the paper we
consider the questions whether, given $\fU$ and $\fV$,
the property $\split(\fU,\fV)$ is preserved under taking finite
or countable unions,
arbitrary subsets, powers or products.
Several interesting problems remain open.
\end{abstract}

\keywords{%
$\gamma$-cover,
$\omega$-cover,
$\tau$-cover,
splitting,
ultrafilter,
$P$-point,
powers,
products,
hereditarity}
\subjclass{03E05, 54D20, 54D80}

\maketitle
%\markboth{\sc Boaz Tsaban}{\sc Selection principles and the minimal tower problem}

%\centerline{\emph{Preliminary version -- comments are welcome}}

\section{Introduction and basic facts}

We consider infinite topological spaces which are homeomorphic to
sets of real numbers (this is the case, e.g., for each separable
and zero-dimensional metric space). We will refer to such spaces
as \emph{sets of reals}. Assume that $X$ is a set of reals. The
following types of ``thick'' covers of $X$ were defined in the
literature and studied under various guises (e.g., \cite{GN, coc1,
coc2, CBC, tau, tautau}). Let $\cU$ be a collection of subsets of
$X$ such that $X$ is not contained in any member of $\cU$. $\cU$
is: \be \i A \emph{large cover} of $X$ if each $x\in X$ is
contained in infinitely many members of $\cU$, \i An
\emph{$\w$-cover} of $X$ if each finite subset of $X$ is contained
in some member of $\cU$, \i A \emph{$\tau$-cover} of $X$ if it is
a large cover of $X$, and for each $x,y\in X$, either $\{U\in\cU :
x\in U, y\nin U\}$ is finite, or $\{U\in\cU : y\in U, x\nin U\}$
is finite; and \i A $\gamma$-cover of $X$ if $\cU$ is infinite,
and each $x\in X$ belongs to all but finitely many members of
$\cU$. \ee Let $\Lambda$, $\Omega$, $\Tau$, and $\Gamma$ denote
the collections of open large covers, $\w$-covers, $\tau$-covers,
and $\gamma$-covers of $X$, respectively. Also, let
$\BL,\BO,\BT,\BG$ (respectively, $\CL,\CO,\CT,\CG$) be the
corresponding \emph{countable Borel} (respectively, \emph{clopen})
covers of $X$. We will informally refer to all these collections
as \emph{collections of thick covers}. It is easy to see that
$$\Gamma\sbst\Tau\sbst\Omega\sbst\Lambda.$$
%and that the analogue assertions for the Borel and clopen cases also
%hold.
Reverse inclusions need not hold. Consider the property
$\binom{\fU}{\fV}$ (read: \emph{$\fU$ choose $\fV$}),
defined for collections of covers $\fU$ and $\fV$,
which asserts that for each cover $\cU\in\fU$
there exists a subcover $\cV\sbst\cU$ such that $\cV\in\fV$.
Then $\binom{\Lambda}{\Omega}$ never holds \cite{coc2, strongdiags},
and there exist sets of reals which do not satisfy
$\binom{\Tau}{\Gamma}$ and $\binom{\Omega}{\Tau}$ \cite{tau, tautau, ShTb768}.

Assume that $\fU$ and $\fV$ are collections of covers of a space $X$.
The following property was introduced in \cite{coc1}.
\bi
\i[$\split(\fU,\fV)$:] Every cover $\cU\in\fU$ can be split
into two disjoint subcovers $\cV$ and $\cW$ which contain elements of $\fV$.
\ei
Several results about these properties (where $\fU,\fV$ are collections of thick
covers) are scattered in the literature. Some of them relate them to classical
properties. For example, it is known that
the Hurewicz property and Rothberger's property both imply $\split(\Lambda,\Lambda)$, and that
the Sakai property (asserting that each finite power of $X$ has Rothberger's property)
implies $\split(\Omega,\Omega)$ \cite{coc1}.
It is also known that if all finite powers of $X$ have the Hurewicz property,
then $X$ satisfies $\split(\Omega,\Omega)$ \cite{coc7}.
By a recent characterization of the Reznichenko (or: weak Fr\'echet-Urysohn) property of $C_p(X)$
in terms of covering properties of $X$ \cite{SakaiNew},
the Reznichenko property for $C_p(X)$ implies that $X$ satisfies $\split(\CO,\CO)$.

Some other works study these properties \emph{per se} \cite{coc2, splitomega}.
As any infinite subset of a $\gamma$-cover is a $\gamma$-cover, we have that any set of reals
satisfies $\split(\Gamma,\Gamma)$ (and therefore $\split(\Gamma,\fV)$ for
all $\fV\spst\Gamma$) \cite{coc1}.
The properties $\split(\Omega,\Omega)$ and $\split(\Lambda,\Lambda)$
are more restrictive \cite{coc2, splitomega}.

\subsection*{Countable subcovers}
%\section{Classification}

It will be more convenient to work with countable covers instead of
covers of arbitrary size.
Each infinite subset of a $\gamma$-cover of a space is
a $\gamma$-cover of the same space. Therefore any
$\gamma$-cover contains a countable $\gamma$-cover.
It is also true (but less trivial) that
every $\omega$-cover of a set of reals $X$
contains a countable $\omega$-cover of $X$ \cite{GN}.

\begin{prop}\label{countablelarge}
Assume that $X$ is a set of reals and $\cU$ is an open large cover of $X$.
Then $\cU$ contains a countable large cover of $X$.
\end{prop}
\begin{proof}
For a cover $\cV$ of a set $Y$ write
$$\cV(Y) = \{y\in Y : y\in V\mbox{ for infinitely many }V\in\cV\}.$$
Write $X_0 = X$. As $X_0$ is Lindel\"of, $\cU$ contains a countable subcover $\cU_0$ of $X_0$.
Set $X_1 = X\sm\cU_0(X_0)$. Then $\cU\sm\cU_0$ is a large cover of $X_1$ (which
is Lindel\"of) and therefore contains a countable subcover $\cU_1$ of $X_1$.
Continue in this manner to define, for each $n$, the sets
$X_n,\cU_n$ such that $X_n=X\sm\cU_{n-1}(X_{n-1})$, and $\cU_n=\cU\sm\Union_{k<n}\cU_k$ is a cover of $X_n$.
Let $X' = \bigcap_nX_n$
and $\cV = \Union_n\cU_n$. As each $\cU_n$ is a countable cover of $X'$ and the sets $\cU_n$, $n\in\N$,
are pairwise disjoint, $\cV$ is a countable large cover of $X'$. For each $x\in X\sm X'$ there exists $n$
such that $x\in\cU_n(X_n)$. Thus $\cV$ is also a large cover of $X\sm X'$, and therefore of $X$.
\end{proof}

We now prove the analogue fact for $\tau$-covers.
\begin{prop}\label{ctbltau}
Assume that $X$ is a set of reals and $\cU$ is an open $\tau$-cover of $X$.
Then $\cU$ contains a countable $\tau$-cover of $X$.
\end{prop}
Proposition \ref{ctbltau} follows from Proposition \ref{countablelarge} and
the following observation, which is of independent importance.
\begin{lem}\label{largeistau}
Assume that $\cU$ is a $\tau$-cover of $X$
and that $\cV\sbst\cU$ is a large cover of $X$.
Then $\cV$ is a $\tau$-cover of $X$.
\end{lem}
\begin{proof}
Assume that $\cU$ is a $\tau$-cover of $X$ and $\cV\sbst\cU$ is a large cover
of $X$. We need only check that for each $x,y\in X$, one of the
sets $\{U\in\cV : x\in U, y\nin U\}$ and
$\{U\in\cV : y\in U, x\nin U\}$ is finite. But
these are subsets of  $\{U\in\cU : x\in U, y\nin U\}$ and
$\{U\in\cU : y\in U, x\nin U\}$, respectively.
\end{proof}

We may therefore assume that all the covers
we consider are countable.
Consequently, the following, where an arrow denotes inclusion, holds:
$$\begin{matrix}
\BG      & \impl & \BT      & \impl & \BO      & \impl & \BL      \\
\uparrow &       & \uparrow &       & \uparrow &       & \uparrow \\
\Gamma   & \impl & \Tau     & \impl & \Omega   & \impl & \Lambda  \\
\uparrow &       & \uparrow &       & \uparrow &       & \uparrow \\
\CG      & \impl & \CT      & \impl & \CO      & \impl & \CL
\end{matrix}$$
As the property $\split(\fU,\fV)$ is monotonic in its first variable and
anti-monotonic in its second variable, we have that for each $x,y\in\{\Gamma, \Tau, \Omega, \Lambda\}$,
$$\split(\B_x,\B_y)\impl\split(x,y)\impl\split(C_x,C_y).$$
Following the mainstream of papers dealing with collections of thick covers, we will
be mostly interested in the splittability properties in the case of (general) open covers,
but we will often use the fact that these properties are ``sandwiched'' between the corresponding
Borel and clopen properties in order to derive theorems about them.

\subsection*{A Ramseyan property}
It is well known \cite{coc1, splitomega} that being an $\omega$-cover is
a Ramsey theoretic property: If an $\omega$-cover is partitioned into
finitely many pieces, then at least one of the pieces is an $\omega$-cover.
The same is true for $\tau$-covers.

\begin{cor}\label{ramseytau}
Assume that $\cU=\cU_1\cup\dots\cup\cU_k$ is a $\tau$-cover of $X$.
Then at least one of the sets $\cU_i$ is a $\tau$-cover of $X$.
\end{cor}
\begin{proof}
$\cU$ is, in particular, an $\omega$-cover of $X$.
Now use the corresponding fact for $\omega$-covers and Lemma \ref{largeistau}.
\end{proof}

An \emph{ultrafilter on $\N$} is a
family $U$ of subsets of $\N$ that
is closed under taking supersets, is closed under
finite intersections, does not contain the empty set as an element,
and for each $a\sbst\N$, either $a\in U$ or $\N\sm a\in U$.
An ultrafilter $U$ on $\N$ is \emph{nonprincipal} if it is not of the
form $\{a\sbst\N : n\in a\}$ for any $n$.

\begin{cor}
Assume that $\cU=\seq{U_n}$ is a $\tau$-cover of a space $X$ which cannot be split
into two $\tau$-covers of $X$. Then
$$U = \{a\sbst\N : \cV=\{U_n\}_{n\in a}\mbox{ is a $\tau$-cover of }X\}$$
is a nonprincipal ultrafilter on $\N$.
\end{cor}
\begin{proof}
This follows from Corollary \ref{ramseytau}, as in \cite{splitomega}.
Alternatively, use Lemma \ref{largeistau}
and the corresponding assertion for $\omega$-covers, which is also true
\cite{splitomega}.
\end{proof}

\part{Classification}

\section{Equivalences and implications}
We begin with the following complete array of properties
(where an arrow denotes implication):
\newcommand{\aru}{\ar[r]\ar[u]}
$$\xymatrix{
\split(\Lambda,\Lambda)\ar[r]&\split(\Omega,\Lambda)\ar[r]&\split(\Tau,\Lambda)\ar[r]&\split(\Gamma,\Lambda)\\
\split(\Lambda,\Omega)\aru&\split(\Omega,\Omega)\aru&\split(\Tau,\Omega)\aru&\split(\Gamma,\Omega)\ar[u]\\
\split(\Lambda,\Tau)\aru&\split(\Omega,\Tau)\aru&\split(\Tau,\Tau)\aru&\split(\Gamma,\Tau)\ar[u]\\
\split(\Lambda,\Gamma)\aru&\split(\Omega,\Gamma)\aru&\split(\Tau,\Gamma)\aru&\split(\Gamma,\Gamma)\ar[u]
}$$

\forget
$$\begin{matrix}
\split(\Lambda, \Lambda) & \to & \split(\Omega, \Lambda) & \to & \split(\Tau, \Lambda) & \to & \split(\Gamma, \Lambda)\\
\uparrow     &     & \uparrow     &     & \uparrow     &     & \uparrow    \\
\split(\Lambda, \Omega) & \to & \split(\Omega, \Omega) & \to & \split(\Tau, \Omega) & \to & \split(\Gamma, \Omega)\\
\uparrow     &     & \uparrow     &     & \uparrow     &     & \uparrow    \\
\split(\Lambda, \Tau) & \to & \split(\Omega, \Tau) & \to & \split(\Tau, \Tau) & \to & \split(\Gamma, \Tau)\\
\uparrow     &     & \uparrow     &     & \uparrow     &     & \uparrow    \\
\split(\Lambda, \Gamma) & \to & \split(\Omega, \Gamma) & \to & \split(\Tau, \Gamma) & \to & \split(\Gamma, \Gamma)
\end{matrix}$$
\forgotten

As we already mentioned in Section 1,
all properties in the last column are trivial in the sense that all sets
of reals satisfy them.
On the other hand, all properties but the top one in the first column imply $\binom{\Lambda}{\Omega}$ and are
therefore trivial in the sense that no infinite set of reals satisfies any of them.

\begin{thm}
The properties
$\split(\Tau , \Tau)$, $\split(\Tau , \Omega)$, and $\split(\Tau , \Lambda)$ are
equivalent.
\end{thm}
\begin{proof}
This is an immediate consequence of Lemma \ref{largeistau}.
\end{proof}

Thus, removing trivialities and equivalences, we are left with the following properties.
$$\xymatrix{
\split(\Lambda,\Lambda)\ar[r]&\split(\Omega,\Lambda)\ar[r]&\split(\Tau,\Tau)\\
&\split(\Omega,\Omega)\ar[u]\\
&\split(\Omega,\Tau)\ar[u]\\
&\split(\Omega,\Gamma)\aru&\split(\Tau,\Gamma)\ar[uuu]
}$$

\forget
$$\begin{matrix}
\split(\Lambda, \Lambda) & \to & \split(\Omega, \Lambda) & \to & \split(\Tau, \Tau)\\
     &     & \uparrow\\
                         &     & \split(\Omega, \Omega)\\
                              &     & \uparrow &                 & \smash{\big{\uparrow}}\\
                           &  & \split(\Omega, \Tau)\\
                             &     & \uparrow\\
                            &  & \split(\Omega, \Gamma) & \to & \split(\Tau,\Gamma)\\
\end{matrix}$$
\forgotten

The following easy cancellation laws can be added to those given in \cite{tautau}.
\begin{prop}\label{cancellation}
If $\fW\sbst\fV\sbst\fU$, then:
\be
\i $\binom{\fU}{\fV}\cap\split(\fV,\fW) = \split(\fU,\fW)$; and
\i $\split(\fU,\fV)\cap\binom{\fV}{\fW} = \split(\fU,\fW)$.
\ee
\end{prop}

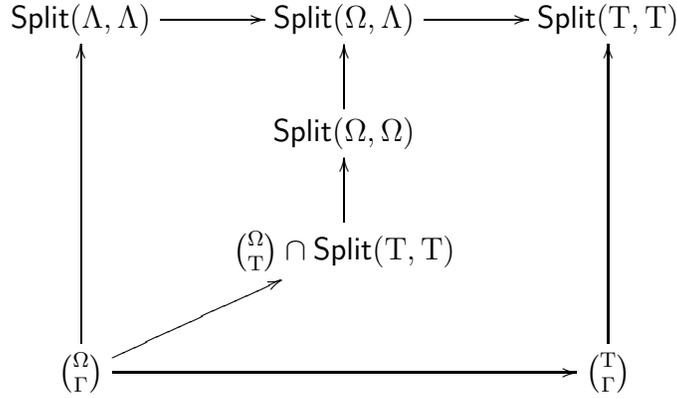
\begin{figure}[!h]
$$\xymatrix{
\split(\Lambda, \Lambda) \ar[r] & \split(\Omega, \Lambda) \ar[r] & \split(\Tau, \Tau)\\
                           & \split(\Omega, \Omega)\ar[u]\\
                          & \binom{\Omega}{\Tau}\cap\split(\Tau, \Tau)\ar[u]\\
\binom{\Omega}{\Gamma} \ar[uuu]\ar[ur]\ar[rr]     & & \binom{\Tau}{\Gamma}\ar[uuu]\\
}$$
\caption{The surviving properties.}\label{survopen}
\end{figure}

\begin{cor}
The following equivalences hold:
\be
\i $\split(\Omega,\Gamma) = \binom{\Omega}{\Gamma}$; and
\i $\split(\Tau,\Gamma) = \binom{\Tau}{\Gamma}$.
\ee
\end{cor}
\begin{proof}
As every set of reals satisfies $\split(\Gamma,\Gamma)$, we have by
Proposition \ref{cancellation} that
$$\binom{\Omega}{\Gamma} = \binom{\Omega}{\Gamma}\cap\split(\Gamma,\Gamma) = \split(\Omega,\Gamma).$$
The proof of the second assertion is similar.
\end{proof}
$\binom{\Omega}{\Gamma}$ is the famous $\gamma$-property introduced by Gerlits and Nagy
in \cite{GN}. The property $\binom{\Tau}{\Gamma}$ was studied in \cite{tautau}.
The property $\split(\Omega, \Tau)$
can also be expressed in terms of other properties:
By Proposition \ref{cancellation},
$$\split(\Omega, \Tau) = \binom{\Omega}{\Tau}\cap\split(\Tau, \Tau).$$
Recall from Section 1 that the Hurewicz property implies $\split(\Lambda,\Lambda)$.
It is well known that the $\gamma$-property implies the Hurewicz property.
Figure \ref{survopen} summarizes our status.
The figures for the clopen and countable Borel cases
are the same, and, as noted before, each property in the Borel case implies
the corresponding property in the open case, which in turn implies
the corresponding property in the clopen case.

\forget
$$\begin{matrix}
\split(\Lambda, \Lambda) & \to & \split(\Omega, \Lambda) & \to & \split(\Tau, \Tau)\\
     &     & \uparrow\\
                         &     & \split(\Omega, \Omega)\\
\smash{\big{\uparrow}}        &     & \uparrow &                 & \smash{\big{\uparrow}}\\
                           &   & \binom{\Omega}{\Tau}\cap\split(\Tau, \Tau)\\
                             & \nearrow    & \\
\binom{\Omega}{\Gamma}        &  & \longrightarrow & & \binom{\Tau}{\Gamma}\\
\end{matrix}$$
\forgotten

\section{Combinatorial characterizations}

In this section we give combinatorial characterizations for all
splitting properties in the cases where the collections of covers
are clopen or countable Borel. These characterizations will be
used in the coming sections to rule out most of the nonexisting
implications between the properties in Figure \ref{survopen}.

We first set the required terminology.
The Cantor space $\Cantor$ of infinite binary sequences is equipped with
the product topology.
Identify $\Cantor$ with $P(\N)$ by characteristic functions.
Then the sets $O_n=\{a\in P(\N) : n\in a\}$ and their complements
form a clopen subbase for the topology of $P(\N)$.
Consider the subspace $\inf$ of $P(\N)$ consisting of the infinite
sets of natural numbers.
For $a,b\in\inf$, we write
$a\as b$ if $a\sm b$ is finite.

A family $Y\sbst\inf$ is \emph{centered} if it is closed under taking
finite intersections.
A family $Y\sbst\inf$ is \emph{reaping} if for each $a\in\inf$ there
exists $y\in Y$ such that $y\as a$ or $y\as\N\sm a$.
Assume that $U$ is a nonprincipal ultrafilter on $\N$.
Observe that $U$ cannot contain a finite set
as an element. Thus, $U$ is a subset of $\inf$.
(Moreover, all cofinite sets belong to $U$
and therefore $U$ is closed under finite
modifications of its elements.)
A family $B\sbst\inf$
is a \emph{base} for $U$ if
$$U = \{a\in\inf : (\exists b\in B)\ b\as a\}.$$
(Consequently,
a family $B\sbst\inf$ is a base for a nonprincipal ultrafilter on $\N$
if, and only if, $B$ is centered and reaping.)
Finally, a family $B\sbst\inf$ is a \emph{subbase} for a nonprincipal ultrafilter $U$ on $\N$ if
$$U = \{a\in\inf : (\exists k)(\exists b_1,\dots,b_k\in B)\ b_1\cap\dots\cap b_k\as a\}.$$

The following combinatorial characterizations are given in \cite{coc2}.
\begin{thm}\label{contimages1}
For a set of reals $X$:
\be
\i $X$ satisfies $\split(\CL ,\CL)$
if, and only if, every continuous image of $X$ in $\inf$ is not a reaping
family.
\i $X$ satisfies $\split(\CO,\CO)$ if, and only if, every
continuous image of $X$ in $\inf$ is not a subbase for a nonprincipal ultrafilter on $\N$.
\ee
\end{thm}

By the same reasoning (see the proof of Theorem \ref{images2} below), one can prove the following.
\begin{thm}\label{borelimages1}
For a set of reals $X$:
\be
\i $X$ satisfies $\split(\BL ,\BL)$
if, and only if, every Borel image of $X$ in $\inf$ is not a reaping
family.
\i $X$ satisfies $\split(\BO,\BO)$ if, and only if, every
Borel image of $X$ in $\inf$ is not a subbase for a nonprincipal ultrafilter on $\N$.
\ee
\end{thm}

\begin{cor}\label{borelopen1}
For a set of reals $X$:
\be
\i $X$ satisfies $\split(\BL ,\BL)$
if, and only if, every Borel image of $X$ satisfies $\split(\CL,\CL)$.
\i $X$ satisfies $\split(\BO,\BO)$ if, and only if, every
Borel image of $X$ satisfies $\split(\CO,\CO)$.
\ee
\end{cor}

We now give combinatorial characterizations for $\split(\CO,\CL)$ and
$\split(\CT,\CT)$. These characterizations as well as the above-mentioned
ones follow from the following lemma.

%\begin{definition}
With each countable cover of $X$ enumerated bijectively as
$\cU=\{U_n\}_{n\in a}$, where $a\sbst\N$, we associate a function
$h_\cU: X\to P(\N)$, defined
by $h_\cU(x) = \{ n\in a : x\in U_n\}$.
%\end{definition}
Note that $h_\cU$ is a Borel function
whenever $\cU$ is a Borel cover of $X$, and $h_\cU$ is continuous
whenever $\cU$ is a clopen cover of $X$.

An element $a\in\inf$ is a \pi{} of a family
$Y\sbst\inf$ if for each $y\in Y$, $a\as y$.
We will need the following minor extension of the corresponding lemma
from \cite{tau}.
\begin{lem}\label{notions}
Assume that $\cU=\{U_n\}_{n\in a}$, where $a\sbst\N$, is a cover of $X$.
\be
\i $\cU$ is a large cover of $X$ if, and only if, $h_\cU[X]\sbst\inf$.
\i $\cU$ is an $\w$-cover of $X$ if, and only if, $h_\cU[X]$ is centered.
\i $\cU$ is a $\tau$-cover of $X$ if, and only if, $h_\cU[X]\sbst\inf$ and is
linearly ordered by $\as$.
\i $\cU$ contains a $\gamma$-cover of $X$ if, and only if, $h_\cU[X]$ has a \pi{}.
\ee
Moreover, if $f:X\to P(\N)$ is any function, and
$\cV=\seq{O_n}$ is the above-mentioned clopen cover of $P(\N)$,
then $f=h_\cU$ for $\cU=\seq{f^{-1}[O_n]}$.
\end{lem}

For a family $Y\sbst\inf$ and an element $a\in\inf$,
the \emph{restriction} of $Y$ to $a$ is the family
$$Y\| a = \{y\cap a : y\in Y\}.$$
If $Y\|a\sbst\inf$, then we say that this restriction
is \emph{large}.
A nonprincipal ultrafilter $U$ on $\N$ is called a \emph{simple $P$-point}
if there exists a base $B$ for $U$ such that $B$ is linearly ordered
by $\as$. We will call such a base a \emph{simple $P$-point base}.

\begin{thm}\label{images2}
For a set of reals $X$:
\be
\i $X$ satisfies $\split(\CO ,\CL)$
if, and only if, every continuous image of $X$ in $\inf$
is not a base for a nonprincipal ultrafilter on $\N$.
\i $X$ satisfies $\split(\BO ,\BL)$
if, and only if, every Borel image of $X$ in $\inf$ is not a base for
a nonprincipal ultrafilter on $\N$.
\i $X$ satisfies $\split(\CT,\CT)$ if, and only if, every
continuous image of $X$ in $\inf$ is not a simple $P$-point base.
\i $X$ satisfies $\split(\BT,\BT)$ if, and only if, every
Borel image of $X$ in $\inf$ is not a simple $P$-point base.
\ee
\end{thm}
\begin{proof}
Observe that for a cover $\cU=\seq{U_n}$ and any subset
$\cV=\{U_n\}_{n\in a}$ of $\cU$,
$$h_\cV[X] = h_\cU[X]\|a.$$
Assume that $\cU$ is a large cover which cannot be split into
two large subcovers. By Lemma \ref{notions} and the
above observation, this means that
$h_\cU[X]\sbst\inf$, and for each subset $\cV=\{U_n\}_{n\in a}$ of $\cU$,
either $h_\cV[X] = h_\cU[X]\|a$ is not large, or
$h_{\cU\sm\cV}[X] = h_\cU[X]\|(\N\sm a)$ is not large.
In the first case there exists $y\in h_\cU[X]$ such that $y\cap a$ is
finite, that is, $y\as\N\sm a$. Similarly, in the second case there exists
$y\in h_\cU[X]$ such that $y\as a$.
In other words, our assumption on $\cU$ is equivalent to
the fact that $h_\cU[X]$ is reaping.

(1) Assume that $X$ does not satisfy $\split(\CO ,\CL)$ and let $\cU$ be a countable
clopen $\w$-cover of $X$ which cannot be split into two large covers of $X$.
Fix some enumeration of $\cU$.
By Lemma \ref{notions}, $h_\cU[X]$, a continuous image of $X$, is centered.
By the above observation, $h_\cU[X]$ is reaping and therefore a
base for a nonprincipal ultrafilter on $\N$.

To prove the remaining implication, assume that $f:X\to\inf$ is a continuous function
such that $Y = f[X]$ is a base for a nonprincipal ultrafilter on $\N$.
%In particular, $\cF$ is centered and
%therefore $\seq{O_n}$ is an $\omega$-cover of $\cF$.
%Thus, $\cU=\seq{f\inv[O_n]}$ is a clopen $\omega$-cover of $X$.
By Lemma \ref{notions},
$\cU=\seq{f\inv[O_n]}$ is a clopen cover of $X$, and
$f=h_\cU$. Thus, $Y=h_\cU[X]$.
As $Y$ is centered, $\cU$ is an $\omega$-cover of $X$.
As $Y$ is reaping, $\cU$ cannot be split into two large covers of $X$.

(2) is similar to (1).

(3) Recall that $\split(\Tau,\Tau)=\split(\Tau,\Lambda)$.

Assume that $\cU=\seq{U_n}$ is a clopen $\tau$-cover of $X$ which cannot be split into
two large covers of $X$.
$Y = h_\cU[X]\sbst\inf$ and is linearly ordered by $\as$.
In particular, $Y$ is centered.
By the arguments of (1), $Y$ is a base for a nonprincipal ultrafilter on $\N$.
As $Y$ is linearly ordered by $\as$, it is a simple $P$-point base.

Now assume that $f:X\to\inf$ is a continuous function
such that $Y = f[X]$ is a simple $P$-point base. In particular, $Y$ is
linearly ordered by $\as$. As in (1), we get that
$\cU=\seq{f\inv[O_n]}$ is a clopen $\tau$-cover of $X$, and,
as $Y$ is reaping, $\cU$ cannot be split into two large covers.

(4) is similar to (3).
\end{proof}

The proofs of Theorem \ref{images2} and the related arguments for $\split(\CL,\CL)$
and $\split(\CO,\CO)$ actually establish the following extension of Lemma \ref{notions}.

\begin{lem}\label{morenotions}
Assume that $\cU=\{U_n\}_{n\in\N}$ is a cover of $X$.
\be
\i
$\cU$ is a large cover of $X$ which cannot be split into two large covers of $X$
if, and only if,
$h_\cU[X]$ is a reaping family.
\i
$\cU$ is an $\omega$-cover of $X$ which cannot be split into two large covers of $X$
if, and only if,
$h_\cU[X]$ is a base for a nonprincipal ultrafilter on $\N$.
\i
$\cU$ is an $\omega$-cover of $X$ which cannot be split into two $\omega$-covers of $X$
if, and only if,
$h_\cU[X]$ is a subbase for a nonprincipal ultrafilter on $\N$.
\i
$\cU$ is an $\tau$-cover of $X$ which cannot be split into two $\tau$-covers of $X$
if, and only if,
$h_\cU[X]$ is a simple $P$-point base.
\ee
\end{lem}

From Theorem \ref{images2} we get the following.

\begin{cor}
For a set of reals $X$:
\be
\i $X$ satisfies $\split(\BO ,\BL)$
if, and only if, every Borel image of $X$ satisfies $\split(\CO ,\CL)$.
\i $X$ satisfies $\split(\BT,\BT)$ if, and only if, every
Borel image of $X$ satisfies $\split(\CT,\CT)$.
\ee
\end{cor}

The properties $\binom{\CO}{\CG}$, $\binom{\CT}{\CG}$, and $\binom{\CO}{\CT}$
(and therefore $\binom{\CO}{\CT}\cap\split(\CT,\CT)$)
also have combinatorial characterizations which follow from Lemma \ref{notions}.
\begin{thm}\label{morechars}
For a set of reals $X$:
\be
\i $X$ satisfies $\binom{\CO}{\CG}$ if, and only if,
each centered continuous image of $X$ in $\inf$ has a \pi{} \cite{RECLAW}.
\i $X$ satisfies $\binom{\CT}{\CG}$ if, and only if,
each $\as$-linearly ordered continuous image of $X$ in $\inf$
has a \pi{} \cite{tau}.
\i $X$ satisfies $\binom{\CO}{\CT}$ if, and only if,
each centered continuous image of $X$ in $\inf$
has a large restriction which is linearly ordered by $\as$ \cite{tautau}.
\ee
\end{thm}
The analogue Borel version of Theorem \ref{morechars} also holds \cite{CBC, tautau}.

\section{Special elements}\label{specialsets}

Sets which are continuous images of Borel sets are called \emph{analytic}.
In \cite{splitomega} it is proved that any analytic set of reals satisfies $\split(\Omega, \Omega)$.
It is well known that analytic sets can also be defined as sets which are
Borel images of the Cantor space $\{0,1\}^\N$.
Consequently, analytic sets are closed under taking Borel images.

\begin{prop}\label{infexample}
~\be
\i Every analytic set of reals satisfies $\split(\BO,\BO)$ as well as
$\binom{\BT}{\BG}$.
\i The analytic set $\inf$ does not satisfy $\split(\CL,\CL)$, and it does not
satisfy $\binom{\CO}{\CT}$ either.
\i $\split(\BO,\BO)\cap\binom{\BT}{\BG}$ does not imply $\split(\CL,\CL)\cup\binom{\CO}{\CT}$.
\ee
\end{prop}
\begin{proof}
(1) Assume that $X$ is an analytic set of reals.
Then each Borel image $Y\sbst\inf$ of $X$ is analytic and therefore satisfies
$\split(\CO,\CO)$. By Corollary \ref{borelopen1}, $X$ satisfies
$\split(\BO,\BO)$. The second assertion was proved in \cite{tautau}.

(2) The first assertion is an immediate consequence of Theorem \ref{contimages1}.
(This is also proved in \cite{splitomega}.)
It remains to prove the second assertion.
It is well known that $\inf$ does not have the $\gamma$-property (which
implies measure zero) \cite{GN}, and that
for separable zero-dimensional metric spaces (this is the case for $\inf$),
$\binom{\Omega}{\Gamma}=\binom{\CO}{\CG}$ (an open $\omega$-cover can
be refined to a clopen $\omega$-cover) \cite{RECLAW}.
Thus $\inf$ does not satisfy $\binom{\CO}{\CG}$.
As $\binom{\CO}{\CT}\cap\binom{\CT}{\CG}=\binom{\CO}{\CG}$,
we have by (1) that $\inf$ does not satisfy $\binom{\CO}{\CT}$.

(3) Follows from (1) and (2).
\end{proof}

Thus, no arrow can be added from
$\split(\Omega,\Omega)$ or from $\binom{\Tau}{\Gamma}$ to any of $\split(\Lambda,\Lambda)$
and $\binom{\Omega}{\Tau}\cap\split(\Tau,\Tau)$.

\begin{cor}\label{LLminusOmOm}
The closed unit interval $I=[0,1]$
satisfies $\split(\Lambda,\Lambda)$,
$\split(\Omega,\Omega)$, and $\binom{\Tau}{\Gamma}$, but
does not satisfy $\binom{\Omega}{\Tau}$.
\end{cor}
\begin{proof}
The Hurewicz property implies $\split(\Lambda,\Lambda)$,
and $\sigma$-compact sets have the Hurewicz property.
Moreover, as $\sigma$-compact sets of reals are $F_\sigma$, they satisfy
$\split(\Omega,\Omega)$ as well as $\binom{\Tau}{\Gamma}$
by Proposition \ref{infexample}.
Finally, the unit interval does not satisfy $\binom{\Omega}{\Gamma}$
and the required assertion follows as in the proof of
Proposition \ref{infexample}.
\end{proof}

In particular, we cannot add an arrow from $\split(\Lambda,\Lambda)$
to $\binom{\Omega}{\Tau}\cap\split(\Tau,\Tau)$ in Figure \ref{survopen}.

One may wonder whether all examples in $\split(\Lambda,\Lambda)\cap\split(\Omega,\Omega)$
are $\sigma$-compact.
The answer for this is negative.

\begin{thm}
There exists a set of reals $X$ such that $X$ is not $\sigma$-compact, and
$X$ satisfies $\split(\Lambda,\Lambda)$ and $\split(\Omega,\Omega)$.
\end{thm}
\begin{proof}
In \cite{ideals} a set of reals $X$ is constructed which is not $\sigma$-compact, and
such that all finite powers of $X$ have the Hurewicz property.
In \cite{coc7} it is proved that any set with this property satisfies
$\split(\Omega,\Omega)$. As $X$ has the Hurewicz property, it also
satisfies $\split(\Lambda,\Lambda)$.
\end{proof}

Corollary \ref{LLminusOmOm} does not rule out the possibility that
$\split(\BO,\BO)$ implies $\binom{\CO}{\CT}$. This nonimplication will
be proved in the next section.

\section{Consistency results}
Thus far we have not used any special hypotheses beyond the usual axioms of
mathematics (ZFC).
In this section we obtain several nonimplications by applying set-theoretic
consistency results.

\begin{thm}\label{noPpts}
It is consistent that all sets of reals satisfy $\split(\BT,\BT)$.
In particular, $\split(\BT,\BT)$ does not imply any of $\split(\CO,\CL)$ and $\binom{\CT}{\CG}$.
\end{thm}
\begin{proof}
In \cite{properforcing} (see also \cite{jubar}) a model of set theory is constructed
where there exist no simple $P$-points.
By Theorem \ref{images2}(4), every set of reals in this model satisfies
$\split(\BT,\BT)$. By Zorn's Lemma
there exists a nonprincipal ultrafilter $U$ on $\N$.
By Theorem \ref{images2}(1), $U$ does not satisfy $\split(\CO,\CL)$.
Also, one can construct by transfinite induction a $\as$-linearly ordered
family $Y\sbst\inf$ which has no \pi{}.
By Theorem \ref{morechars}(2), $Y$ does not satisfy $\binom{\CT}{\CG}$.
\end{proof}

A natural question is whether $\split(\Tau,\Tau)$ is, like $\split(\Gamma,\Gamma)$,
trivial in the sense that all sets of reals satisfy this property.
It is easy to construct, assuming \CH{} (or just $\t=\c$ -- see definitions below),
a $\as$-decreasing sequence $\<a_\alpha : \alpha<\c\>$ such that for each
For each $a\sbst\N$, there exists
$\alpha$ such that either $a_\alpha\as a$ or $a_\alpha\as \N\sm a$ \cite{tau}.
Clearly such a sequence forms a simple $P$-point base, and,
by Theorem \ref{images2}, does not satisfy $\split(\Tau,\Tau)$.
The following shows a bit more than that (at the cost of using a very deep result).
Let $\c$ denote the cardinality
of the continuum. In \cite{ShBl} a model of set theory is constructed in which
$\c=\aleph_2$ and there exist two simple $P$-points with bases of cardinalities
$\aleph_1$ and $\aleph_2$.
\begin{cor}\label{ShBlmodel}
It is consistent that $\c=\aleph_2$ and there exist sets of reals $X$ and $Y$
of cardinalities $\aleph_1$ and $\aleph_2$, respectively, which do not
satisfy $\split(\Tau,\Tau)$.
\end{cor}

In order to proceed, we introduce several cardinal characteristics of the
continuum and some of their properties (see \cite{vD, Blass} for details and proofs).
Let $\r$ denote the minimal cardinality of a reaping family, and
$\u$ denote the minimal cardinality of a base for a nonprincipal ultrafilter
on $\N$.
Then $\r\le\u$.
The \emph{critical cardinality} of a property $\bP$ of sets of reals, $\non(\bP)$,
is the minimal cardinality of a set of reals which does not satisfy this property.
In \cite{coc2} it is deduced from Theorem \ref{contimages1} that
$\non(\split(\Lambda,\Lambda))=\r$, and $\non(\split(\Omega,\Omega))=\u$.
(These results also hold in the clopen and Borel cases.)
By Theorem \ref{images2}, we have the following.
\begin{thm}
The critical cardinalities of the classes
$\split(\BO,\BL)$, $\split(\Omega,\Lambda)$, and $\split(\CO,\CL)$ are all equal to $\u$.
\end{thm}

Let $\p$ denote the minimal cardinality of a centered family in $\inf$ which
does not have a \pi{}.
In \cite{RECLAW, CBC, tautau} it is shown that the critical cardinalities of
$\binom{\BO}{\BG}$, $\binom{\Omega}{\Gamma}$, $\binom{\CO}{\CG}$,
$\binom{\BO}{\BT}$, $\binom{\Omega}{\Tau}$, and $\binom{\CO}{\CT}$
are all equal to $\p$.
\begin{cor}
The critical cardinalities of
$\binom{\BO}{\BT}\cap\split(\BT,\BT)$,
$\binom{\Omega}{\Tau}\cap\split(\Tau,\Tau)$, and
$\binom{\CO}{\CT}\cap\split(\CT,\CT)$
are all equal to $\p$.
\end{cor}
\begin{proof}
All these properties are implied by
$\binom{\BO}{\BG}$ (whose critical cardinality is $\p$), and imply
$\binom{\CO}{\CT}$ (whose critical cardinality is also $\p$).
\end{proof}
A \emph{tower} of length $\kappa$ is a $\as$-decreasing
sequence $\<a_\alpha : \alpha<\kappa\>$ of elements of $\inf$,
which has no \pi{}. Let $\t$ denote the minimal cardinality of
a tower. In \cite{tau, tautau} it is deduced from Theorem \ref{morechars}
and its Borel version that
%By Theorem \ref{morchars}(2), no tower
%satisfies $\binom{\CT}{\CG}$ \cite{tau}.
the critical cardinalities of the classes
$\binom{\BT}{\BG}$, $\binom{\Tau}{\Gamma}$, and $\binom{\CT}{\CG}$ are equal to $\t$.
%It is easy to see that $\aleph_1\le\p\le\t$.
The following diagram summarizes the critical cardinalities of the properties we study
(observe that by Theorem \ref{noPpts}, the critical cardinality of $\split(\Tau,\Tau)$
is undefined).
$$
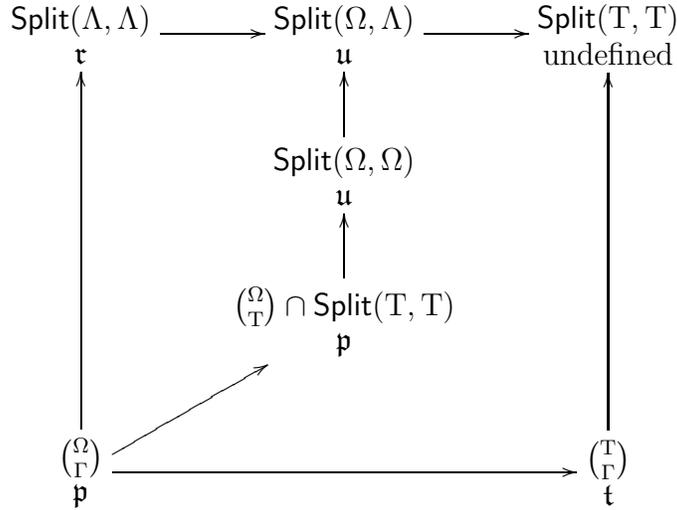
\xymatrix{
\sr{\split(\Lambda, \Lambda)}{\r} \ar[r] & \sr{\split(\Omega, \Lambda)}{\u} \ar[r] & \sr{\split(\Tau, \Tau)}{\mbox{undefined}}\\
                           & \sr{\split(\Omega, \Omega)}{\u}\ar[u]\\
                          & \sr{\binom{\Omega}{\Tau}\cap\split(\Tau, \Tau)}{\p}\ar[u]\\
\sr{\binom{\Omega}{\Gamma}}{\p} \ar[uuu]\ar[ur]\ar[rr]     & & \sr{\binom{\Tau}{\Gamma}}{\t}\ar[uuu]\\
}$$

Let $\h$ denote the \emph{distributivity number}.
For our purposes the definition of $\h$ is not important;
we need only quote the result that $\h\le\r$.
The following theorem strengthens Theorem \ref{noPpts}.

\begin{thm}
There exists a single model of set theory that witnesses the following
facts:
\be
\i $\split(\BT,\BT)$ does not imply any of $\split(\CO,\CL)$ and $\binom{\CT}{\CG}$; and
\i $\split(\BL,\BL)\cap\split(\BO,\BO)$ does not imply any of
$\binom{\CO}{\CT}$ and $\binom{\CT}{\CG}$.
\ee
\end{thm}
\begin{proof}
In \cite{dordal} a model of set theory is constructed in which
$\h = \c = \aleph_2$ but there are
no towers of length $\aleph_2$.
As $\h\le\r$, $\r=\u=\c=\aleph_2$ in this model.
\begin{lem}
There exist no simple $P$-points in this model.
\end{lem}
\begin{proof}
Assume that $B\sbst\inf$ is a simple $P$-point base.
Then $|B|\ge\u$. As $\u=\c=\aleph_2$,
$|B|=\aleph_2$, and a cofinal $\as$-decreasing subset
of $\cF$ would be a tower of length $\aleph_2$,
a contradiction.
\end{proof}
Thus, in this model all sets of reals satisfy $\split(\BT,\BT)$.

As there are no towers of length $\aleph_2$ in this model, we have that $\p=\t=\aleph_1$.
Thus there exist sets of reals $X$ and $Y$  of cardinality $\aleph_1$ which do not
satisfy $\binom{\CO}{\CT}$ and $\binom{\CT}{\CG}$, respectively.
As $\aleph_1<\r\le\u$, $X$ and $Y$ satisfy $\split(\BL,\BL)$ as well as $\split(\BO,\BO)$.
\end{proof}

We now prove that $\split(\Lambda,\Lambda)$ does not imply $\split(\Omega,\Omega)$.
The \emph{additivity number} of a collection (or a property) $\cI$ of sets of reals
is
$$\add(\cI) = \min\{|\cF| : \cF\sbst\cI\mbox{ and }\cup\cF\nin\cI\},$$
and the \emph{covering number} of $\cI$ is
$$\cov(\cI) = \min\{|\cF| : \cF\sbst\cI\mbox{ and }\cup\cF=\R\}.$$
Let $\M$ denote the collection of meager (i.e., first category) sets of real numbers.
By the Baire's category theorem, $\add(\M)\le\cov(\M)$.
Assume that $\kappa$ is an uncountable cardinal. A set of reals $L$ is a \emph{$\kappa$-Luzin set}
if $|L|\ge\kappa$ and for each meager set $M$, $|L\cap M|<\kappa$.

\begin{thm}\label{specialLuzin}
Assume \CH{} (or just $\add(\M)=\c$).
Then there exists an $\add(\M)$-Luzin set $L$ which satisfies
$\split(\BL,\BL)$ but not $\split(\CO,\CO)$.
\end{thm}
\begin{proof}
In \cite{CBC} it is proved that if $L$ is an $\add(\M)$-Luzin set,
then each Borel image of $L$ satisfies Rothberger's property.
As Rothberger's property implies $\split(\CL,\CL)$ \cite{coc1},
we have by Corollary \ref{borelopen1} that $L$ satisfies $\split(\BL,\BL)$.

It therefore suffices to construct a $\add(\M)$-Luzin set which is
a subbase for a nonprincipal ultrafilter on $\N$.
To this end, fix a nonprincipal ultrafilter $U$ on $\N$.
It is well known that nonprincipal ultrafilters on $\N$
do not have the Baire property, and in particular are nonmeager \cite{jubar}.
It is therefore conceivable that the following holds.
\begin{lem}\label{UminusM}
Assume that $U$ is a nonprincipal ultrafilter on $\N$ and that
$M\sbst\inf$ is meager. Then $U\sm M$ is a subbase for $U$.
In fact, for each $a\in U$ there exist $a_0,a_1\in U\sm M$ such
that $a_0\cap a_1\sbst a$.
\end{lem}
Replying to a question of ours, Shelah gave a proof for this lemma.
To simplify the proof, we make some translation.
Recall that $\inf$ is a subspace of $P(\N)$ whose topology is
defined by its identification with $\Cantor$.
It is well known \cite{jubar, Blass} that for each meager subset $M$ of $\Cantor$
there exist $x\in\Cantor$ and a strictly increasing function
$f\in\NN$ such that
$$M\sbst\{y\in\Cantor : (\A^\oo n)\ y\|\left [f(n),f(n+1)\right )\neq x\|\left [f(n),f(n+1)\right )
\}$$
where $\A^\oo n$ means ``for all but finitely many $n$''.
%(The set on the right hand side is also meager.)
Translating this to the language of $\inf$, we get that
%for each meager $M\sbst\inf$ there exists a strictly increasing function
%$f\in\NN$, and
for each $n$ there exist disjoint sets $I^n_0$ and $I^n_1$ satisfying
$I^n_0\cup I^n_1 = \left [f(n),f(n+1)\right )$,
such that
\begin{equation}\label{Mcof}
M\sbst\{y\in\inf : (\A^\oo n)\ y\cap I^n_0\neq\emptyset\mbox{ or } I^n_1\not\sbst y\}.
\end{equation}
\begin{proof}[Proof of Lemma \ref{UminusM}]
Assume that the sets $I^n_0,I^n_1$, $n\in\N$, are chosen as in
(\ref{Mcof}).
Let $a$ be an infinite co-infinite subset of $\N$.
Then either $x=\Union_{n\in a}\left [f(n),f(n+1)\right )\nin U$, or else
$x=\Union_{n\in \N\sm a}\left [f(n),f(n+1)\right )\nin U$.
We may assume that the former case holds.
Split $a$ into two disjoint infinite sets $a_1$ and $a_2$.
Then $x_i=\Union_{n\in a_i}\left [f(n),f(n+1)\right )\nin U$ ($i=1,2$).

Assume that $b\in U$. Then $\tilde b = b\sm x = b\cap(\N\sm x)\in U$.
Define sets $y_1,y_2\in U\sm M$ as follows.
\begin{eqnarray*}
%Y_1 & = & \tilde B \cup \Union_{n\in A_1}I^n_0 \cup \Union_{n\in A_2}I^n_1\\
%Y_2 & = & \tilde B \cup \Union_{n\in A_1}I^n_1 \cup \Union_{n\in A_2}I^n_0
y_1 & = & \tilde b \cup \Union_{n\in a_2}I^n_1\\
y_2 & = & \tilde b \cup \Union_{n\in a_1}I^n_1
\end{eqnarray*}
By (\ref{Mcof}), $y_1,y_2\nin M$.
As $y_1,y_2\spst\tilde b$, $y_1,y_2\in U$.
Now, $y_1\cap y_2=\tilde b\sbst b$.
\end{proof}
We now construct the Luzin set $L$.
Enumerate $U$ as $\{a_\alpha : \alpha<\c\}$, and let
$\{M_\alpha : \alpha<\c\}$ be a cofinal family of meager sets in $\inf$
(e.g., the $F_\sigma$ meager sets).
For each $\alpha<\c$ use Lemma \ref{UminusM} to choose
$$a^0_\alpha,a^1_\alpha\in U\sm\Union_{\beta<\alpha}M_\beta$$
such that $a^0_\alpha\cap a^1_\alpha\sbst a_\alpha$.
Then $L=\{a^0_\alpha, a^1_\alpha : \alpha<\c\}$ is as required.
\end{proof}

It is an open problem whether $\binom{\Omega}{\Tau}=\binom{\Omega}{\Gamma}$ \cite{tautau}.
Observe that if $\binom{\Omega}{\Tau}$ implies $\binom{\Tau}{\Gamma}$, then
$\binom{\Omega}{\Tau}=\binom{\Omega}{\Gamma}$. The only remaining classification
problems are stated in the following problem.

\begin{prob}\label{survprobs}
Is the dotted implication (1) (and therefore (2) and (3)) in the following diagram true?
If not, then is the dotted implication (3) true?
$$
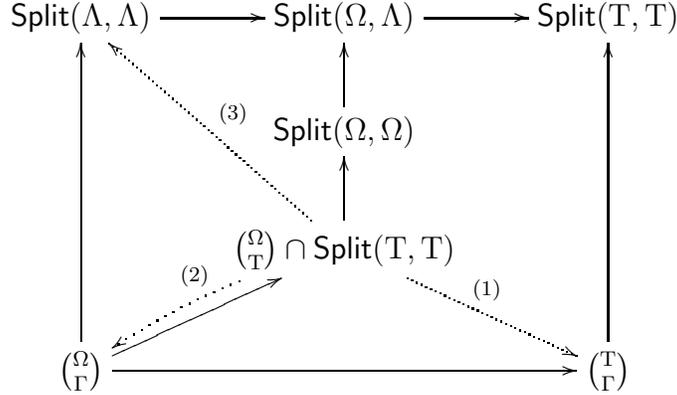
\xymatrix{
\split(\Lambda, \Lambda) \ar[r] & \split(\Omega, \Lambda) \ar[r] & \split(\Tau, \Tau)\\
                           & \split(\Omega, \Omega)\ar[u]\\
& \binom{\Omega}{\Tau}\cap\split(\Tau, \Tau)\ar[u]\ar@{.>}[dr]^{(1)}\ar@/_/@{.>}[dl]_{(2)}\ar@{.>}[uul]_{(3)}\\
\binom{\Omega}{\Gamma} \ar[uuu]\ar[ur]\ar[rr]     & & \binom{\Tau}{\Gamma}\ar[uuu]\\
}$$
\end{prob}

Observe that with regards to the properties $\split(\Lambda, \Lambda)$, $\split(\Omega, \Lambda)$,
$\split(\Tau, \Tau)$, and $\split(\Omega, \Omega)$, the classification is complete.

\part{Preservation of properties}

\section{Unions}\label{unions}

The proof of Theorem \ref{specialLuzin} can be extended to obtain more.
For the proof, we need some notation and results from \cite{huremen2, AddQuad}.
A cover $\cU$ of $X$ is \emph{$\w$-fat} if for each finite $F\sbst X$
and finite family $\cF$ of nonempty open sets, there exists $U\in\cU$ such that
$F\sbst U$, and for each $G\in\cF$ $U\cap G$ is not meager.
In this case, for each finite $F\sbst X$ and finite family $\cF$ of nonempty basic open sets,
the set $\cup\{U\in\cU : F\sbst U\mbox{ and }(\forall O\in\cF)\ U\cap O\nin\M\}$ is comeager, and
for each element $x$ in the intersection of all sets of this form,
$\cU$ is an $\w$-fat cover of $X\cup\{x\}$.
Let $\BOfat$ denote the collection of countable Borel
$\w$-fat covers of $X$.
The following property, which generalizes several
classical properties, was introduced in \cite{coc1}.
\begin{itemize}
\item[$\sone(\fU,\fV)$:]
For each sequence $\seq{\cU_n}$ of members of $\fU$,
there is a sequence
$\seq{U_n}$ such that for each $n$ $U_n\in\cU_n$, and $\seq{U_n}\in\fV$.
\end{itemize}
Then $\non(\sone(\BOfat,\BOfat))\ge\cov(\M)\ge\add(\M)$.
Finally, if $L$ is a $\kappa$-Lusin set such that for each nonempty basic
open set $G$, $|L\cap G|=\kappa$, then every  countable Borel $\w$-cover $\cU$ of $L$
is an $\w$-fat cover of $L$.

\begin{lem}\label{soneBOBO}
$\sone(\BO,\BO)$ implies $\split(\BO,\BO)$ as well as $\split(\BL,\BL)$.
\end{lem}
\begin{proof}
$\sone(\BO,\BO)$, which is closed under taking Borel images,
implies the Sakai property, which implies
$\split(\CO,\CO)$ as well as $\split(\CL,\CL)$.
The assertion follows from Corollary \ref{borelopen1}.
\end{proof}

Observe that a union of two $\add(\M)$-Luzin sets is
again an $\add(\M)$-Luzin set, and therefore satisfies
$\split(\BL,\BL)$. Thus, the following theorem,
apart from showing that the properties $\split(\BO,\BO)$, $\split(\Omega,\Omega)$,
and $\split(\CO,\CO)$ are not additive, also extends Theorem \ref{specialLuzin}.
\begin{thm}\label{notadd}
Assume \CH{} (or just $\add(\M)=\c$). Then there exist two $\add(\M)$-Luzin sets
$L_0$ and $L_1$ satisfying $\sone(\BO,\BO)$ (and therefore $\split(\BO,\BO)$ and $\split(\BL,\BL)$),
such that $L=L_0\cup L_1$ (which satisfies
$\split(\BL,\BL)$) does not satisfy $\split(\CO,\CO)$.
\end{thm}
\begin{proof}
We follow the footsteps of the proof given in \cite{AddQuad}.
Let $U=\{a_\alpha : \alpha<\c\}$ be a nonprincipal ultrafilter on $\N$.
Let $\{M_{\alpha}:\alpha<\c\}$ enumerate all $F_\sigma$ meager sets in $\inf$,
and $\{\seq{\cU^\alpha_n} : \alpha<\c\}$
enumerate all countable sequences of countable families of Borel sets in $\inf$.
Let $\{O_i : i\in\N\}$ and $\{\cF_i : i\in\N\}$ enumerate all nonempty basic open sets
and finite families of nonempty basic open sets, respectively, in $\inf$.

We construct $L_i=\{a^i_\beta : \beta<\c\}$, $i=1,2$,
by induction on $\alpha<\c$ as follows.
At stage $\alpha\ge 0$ set
$X^i_\alpha  = \{a^i_\beta : \beta<\alpha\}$
and consider the sequence $\seq{\cU^\alpha_n}$.
Say that $\alpha$ is $i$-good if for each $n$
$\cU^\alpha_n$ is an $\w$-fat cover of $X^i_\alpha$.
In this case,
by the above remarks there exist elements
$U^{\alpha,i}_n\in\cU^\alpha_n$ such that $\seq{U^{\alpha,i}_n}$ is
an $\w$-fat cover of $X^i_\alpha$.
We make the inductive hypothesis that
for each $i$-good $\beta<\alpha$,
$\seq{U^{\beta,i}_n}$ is an $\w$-fat cover of $X^i_\alpha$.
For each finite $F\sbst X^i_\alpha$, $i$-good $\beta\le\alpha$,
and $m$ define
$$G_i(F,\beta,m)=\cup\{U^{\beta,i}_n : F\sbst U^{\beta,i}_n\mbox{ and }
(\forall O\in\cF_m)\ U^{\beta,i}_n\cap O\nin\M\}.$$
By the inductive hypothesis, $G_i(F,\beta,m)$ is comeager.
Set
$$Y_\alpha=\Union_{\beta<\alpha}M_\beta\cup
\Union_{\begin{matrix}
i<2,\mbox{ $i$-good }\beta\le\alpha\\
m\in\N,\mbox{ Finite }F\sbst X^i_\alpha
\end{matrix}}
(\inf\sm G_i(F,\beta,m)),$$
and $Y_\alpha^* = \{x\in\inf : (\E y\in Y_\alpha)\ x=^* y\}$
(where $x =^* y$ means that $x\as y$ and $y\as x$.)
Then $Y_\alpha^*$ is a union of less than $\add(\M)$ many meager sets,
and is therefore meager.
Use Lemma \ref{UminusM}
to pick $a^0_\alpha,a^1_\alpha\in U\sm Y_\alpha^*$ such that
$a^0_\alpha\cap a^1_\alpha \as a_\alpha$.
Let $k = \alpha \bmod \omega$, and change finitely many elements
of $a^0_\alpha$ and $a^1_\alpha$ so that
they both become members of $O_k$.
Then $a^0_\alpha,a^1_\alpha\in (U\cap O_k)\sm Y_\alpha$,
and $a^0_\alpha \cap a^1_\alpha \as a_\alpha$.
Observe that, by the remarks in the beginning of this section,
the inductive hypothesis remains true for $\alpha$.
This completes the construction.

Clearly $L_0$ and $L_1$ are Luzin sets
and $L_0\cup L_1$ is a subbase for $U$.
We made sure that for each nonempty basic open set $G$,
$|L_0\cap G|=|L_1\cap G|=\c$,
thus $\BO=\BOfat$ for $L_0$ and $L_1$.
By the construction, $L_0,L_1\in\sone(\BOfat,\BOfat)$.
\end{proof}

The properties $\binom{\BT}{\BG}$, $\binom{\Tau}{\Gamma}$, and $\binom{\CT}{\CG}$
 are $\sigma$-additive (their additivity number is exactly $\t$) \cite{tau, tautau}.

We will show that no property between $\binom{\BO}{\BG}$ and $\binom{\CO}{\CT}$
is provably additive.
Let $\bP$ be a property of sets of reals. We say that a set of reals $X$ is
\emph{hereditarily-$\bP$} if all subsets of $X$ satisfy the property $\bP$.

\begin{thm}\label{AunionB}
Assume \CH{}.
There exist disjoint, zero-dimensional sets of reals $A$ and $B$ satisfying $\binom{\BO}{\BG}$,
such that $A\cup B$ does not satisfy $\binom{\CO}{\CT}$.
\end{thm}
\begin{proof}
In \cite{tautau} it is shown that assuming \CH{}, there exist
disjoint, zero-dimensional sets of reals $A\sbst(0,1)$ and $B\sbst(1,2)$ satisfying $\binom{\BO}{\BG}$,
such that $A\cup B$ does not satisfy $\binom{\Omega}{\Tau}$.
In particular, $A\cup B$ does not satisfy $\binom{\Omega}{\Gamma}$.
As $A\sbst(0,1)$ and $B\sbst(1,2)$, $A\cup B$ is zero-dimensional too,
and therefore
%The proofs given in \cite{GM, tau, tautau} essentially prove this theorem.
%Assume \CH{}.
%Then there exists a hereditarily-$\binom{\BO}{\BG}$
%set of reals $X$ of cardinality $\c$ \cite{JORG}.
%As $\binom{\BO}{\BG}$ is closed under taking Borel (continuous is enough) images,
%we may assume that $X\sbst (0,1)$.
%Moreover, it is well known that $\binom{\Omega}{\Gamma}$ implies zero-dimensionality;
%thus $X$ is zero-dimensional.
%
%For $Y\sbst (0, 1)$, write $Y+1 = \{y+1 : y\in Y\}$ for the translation
%of $Y$ by $1$.
%As $|X|=\c$ and only $\c$ many out of the $2^\c$ many subsets of $ X$
%are Borel, there exists a subset $Y$ of $ X$ which
%is not $F_\sigma$ neither $G_\delta$.
%By a theorem of Galvin and Miller \cite{GM},
%for such a subset $Y$ the set $( X\sm Y)\cup (Y+1)$
%does not satisfy $\binom{\Omega}{\Gamma}$.
%As $X$ is zero-dimensional, so is $( X\sm Y)\cup (Y+1)$, and therefore
%this set does not satisfy $\binom{\CO}{\CG}$.
%
%Set $A= X\sm Y$ and $B=Y+1$.
%Then $A$ and $B$ are disjoint, they satisfy
%$\binom{\BO}{\BG}$, and
$A\cup B$ does not satisfy
$\binom{\CO}{\CG}=\binom{\CO}{\CT}\cap\binom{\CT}{\CG}$.
As $\binom{\CT}{\CG}$ is additive,
$A\cup B$ satisfies $\binom{\CT}{\CG}$. Thus,
$A\cup B$ does not satisfy $\binom{\CO}{\CT}$.
\end{proof}

\begin{thm}
The properties $\split(\BT,\BT)$, $\split(\Tau,\Tau)$, and $\split(\CT,\CT)$
are $\sigma$-additive. In fact, they are closed under taking unions of
size less than $\u$.
\end{thm}
This theorem follows Theorem \ref{images2} and the following Ramseyan property.
\begin{lem}
Assume that $\lambda<\u$ and $B = \Union_{\alpha<\lambda}B_\alpha$
is a simple $P$-point base. Then there exists $\alpha<\lambda$
such that $B_\alpha$ is a simple $P$-point base.
\end{lem}
\begin{proof}
Assume that $B$ is a simple $P$-point base
and $U$ is the simple $P$-point it generates.
In particular, $B$ is linearly ordered by $\as$.
We will show that some $B_\alpha$ is a base for $U$.
Assume otherwise.
For each $\alpha<\lambda$ choose $a_\alpha\in U$ that
witnesses that $B_\alpha$ is not a base for $U$,
and $\tilde a_\alpha\in B$ such that $\tilde a_\alpha\as a_\alpha$.
As $B$ is linearly ordered by $\as$,
$\tilde a_\alpha$ is a \pi{} of $B_\alpha$.

The cardinality of the linearly ordered set $Y = \{\tilde a_\alpha : \alpha<\lambda\}$
is smaller than $\u$. Thus it is not a base for $U$ and we can find again
an element $a\in\cF$ which is a \pi{} of $Y$, and therefore of $B$;
a contradiction.
\end{proof}

Using similar ideas, one can prove that the properties in the forthcoming
Theorem \ref{add3} are (finitely) additive. The referee has pointed out to
us that in fact, these properties are $\sigma$-additive. The proof is
almost verbatim the one given by the referee.
\begin{thm}\label{add3}
The properties $\split(\BO,\BL)$, $\split(\Omega,\Lambda)$, and
$\split(\CO,\CL)$ are $\sigma$-additive.
\end{thm}
\begin{proof}
We will prove the open case. The other cases are similar.
\begin{lem}
Assume that $\cU$ is a countable open $\omega$-cover of $Y$ and
that $X\subseteq Y$ satisfies $\split(\Omega,\Lambda)$.
Then $\cU$ can be partitioned into two pieces $\cV$ and $\cW$
such that that $\cW$ is an $\omega$-cover of $Y$ and
each element of $X$ is contained in infinitely many members of
$\cV$.\footnote{Due to our technical requirement in the introduction that $X$ is not
contained in any member of the cover, this does not imply that $\cV$ is a
large cover of $X$.
}
\end{lem}
\begin{proof}
First assume that there does not exist $U\in \cU$ with $X\subseteq
U$.  Then $\cU$ in an $\omega$-cover of $X$. By the splitting
property we can divide it into two pieces each a large cover of
$X$.  Since $\cU$ is an $\omega$-cover of $Y$, one of the pieces
is an $\omega$-cover of $Y$ (see introduction), and the lemma is
proved. If there are only finitely many $U\in \cU$ with
$X\subseteq U$, then $\tilde \cU=\cU\sm\{U\in \cU : X\subseteq
U\}$ is still an $\omega$-cover of $Y$ and we can apply to it the
above argument.

Thus, assume that there are infinitely many $U\in \cU$ with $X\subseteq U$.
Then take a partition of $\cU$ into two pieces such that each piece
contains infinitely many sets $U$ with $X\sbst U$.
One of the pieces must be an $\omega$-cover of $Y$.
\end{proof}
Assume that $Y=\Union_{n\in\N}X_n$ where each $X_n$ satisfies $\split(\Omega,\Lambda)$,
and let $\cU_0$ be an open $\omega$-cover of $Y$.
Given $\cU_n$ an open $\omega$-cover of $Y$, apply the lemma twice to get a partition
$\cU_n=\cV_n^0\cup\cV_n^1\cup\cU_{n+1}$ such that
$\cU_{n+1}$ is an open $\omega$-cover of $Y$ and for each $i=0,1$, each element of
$X_n$ is contained in infinitely many $V\in \cV_n^i$.
Then the families $\cV^i=\Union_{n\in\N}\cV^i_n$, $i=0,1$, are disjoint large covers of $Y$
which are subcovers of $\cU_0$.
\end{proof}

One additivity problem remains open.
\begin{prob}\label{SpLamLamAdd}
Is $\split(\Lambda,\Lambda)$ additive?
\end{prob}

\section{Hereditarity}

We have, implicitly and explicitly, used the following fact in the
preceding sections.
\begin{prop}\label{hered}
For each $x,y\in\{\Lambda,\Omega,\Tau,\Gamma\}$:
\be
\i $\split(C_x,C_y)$ is closed under taking clopen subsets
and continuous images,
\i $\split(x,y)$ is closed under taking closed subsets
and continuous images; and
\i $\split(\B_x,\B_y)$ is closed under taking Borel subsets
and continuous images.
\ee
\end{prop}
\begin{proof}
The proofs for these assertions are standard, see \cite{coc2, CBC}.
\end{proof}

A class $\cI$ of sets of reals is \emph{hereditary} if it is closed
under taking subsets.

\begin{thm}\label{nonhered}
Assume \CH{} (or just $\p=\c$).
Then there exists a set of reals $X$ (of size $\c$)
and a countable subset $Q$ of $X$
such that $X$ satisfies $\binom{\Omega}{\Gamma}$ and $X\sm Q$ does not
satisfy $\split(\CT,\CT)$.
\end{thm}
\begin{proof}
In \cite{ideals}, a subset $X$ of $P(\N)$ is constructed, such
that:
\be
\i $X$ satisfies $\binom{\Omega}{\Gamma}$,
\i $X=P\cup Q$ where $P\sbst\inf$ is linearly ordered by $\as$ and $Q$ is countable; and
\i For each infinite coinfinite subset $a$ of $\N$, there exists
   $x\in P$ such that either $x\as a$, or else $x\as \N\sm a$.
\ee
Consequently, $X\sm Q = P$ is a simple $P$-point base, which, by Theorem \ref{images2},
does not satisfy $\split(\CT,\CT)$.
\end{proof}

\begin{cor}
None of the splittability properties in the open (or clopen) case
implies any of the splittability properties in the Borel case.
\end{cor}
\begin{proof}
Consider the set $X$ given in Theorem \ref{nonhered}.
$X$ satisfies $\binom{\Omega}{\Gamma}$, and
as $Q$ is countable, $X\sm Q$ is a Borel subset of $X$.
By Proposition \ref{hered}, if $X$ satisfied
$\split(\BT,\BT)$, so would $X\sm Q$. In particular,
we would have that $X\sm Q$ satisfies $\split(\CT,\CT)$, a contradiction.
\end{proof}

Despite the above, some classes in the Borel case \emph{are} provably hereditary.
\begin{thm}
$\split(\BL,\BL)$ is hereditary.
\end{thm}
\begin{proof}
This follows from Theorem \ref{borelimages1} and
the fact that each Borel function
defined on a set of reals can be extended to a Borel
function on $\R$ \cite{KURA}.
A direct proof for this is as follows:
Assume that $X$ satisfies $\split(\BL,\BL)$ and that $Y$ is a subset of $X$.
Assume that $\cU$ is a countable Borel cover of $Y$.
Then
$$\cV = \{U\cup (X\sm\cup\cU) : U\in\cU\}$$
is a countable Borel large cover of $X$, and therefore can
be split into two disjoint large subcovers $\cV_1$ and $\cV_2$.
Then $\{V\cap Y : V\in\cV_1\}\sm\{\emptyset\}$ and $\{V\cap Y : V\in\cV_2\}\sm\{\emptyset\}$
are disjoint subsets of $\cU$ and are large covers of $Y$.
\end{proof}

Recently, Miller proved that
no class between $\binom{\BO}{\BG}$ and $\binom{\BO}{\BT}$ is provably
hereditary \cite{MilNonGamma}. In particular, $\binom{\BO}{\BT}\cap\split(\BT,\BT)$
is not provably hereditary.
\begin{prob}\label{BorelHered}
Is any of the remaining classes (namely, $\split(\BO,\BL)$, $\split(\BO,\BO)$, $\split(\BT,\BT)$,
and $\binom{\BT}{\BG}$) provably hereditary?
\end{prob}

\section{Finite powers and products}

The $\gamma$-property $\binom{\Omega}{\Gamma}$ is provably closed under taking
finite powers, but not under taking finite products \cite{coc2}.
This assertion can be extended.
\begin{thm}
No class between $\binom{\BO}{\BG}$ and $\binom{\CO}{\CT}$ is provably
closed under taking finite products.
\end{thm}
\begin{proof}
The proof for this is as in \cite{GM}.
Assume \CH{}, and let $A$ and $B$ be as in Lemma \ref{AunionB}.
Assume that $A\x B$ satisfies $\binom{\CO}{\CT}$.
Fix $a\in A$ and $b\in B$.
As $A$ and $B$ are zero-dimensional,
The set $X = (A\x\{b\})\cup(\{a\}\x B)$ is a clopen subset of $A\x B$
and therefore satisfies $\binom{\CO}{\CT}$ too.
But as $A$ and $B$ are disjoint,
this set is homeomorphic to $A\cup B$, which does not satisfy
$\binom{\CO}{\CT}$, a contradiction.
\end{proof}

In particular, $\binom{\Omega}{\Tau}\cap\split(\Tau,\Tau)$ is not
provably closed under taking finite products.
We do not know whether this property is provably closed under
taking finite \emph{powers}. In fact, we cannot even answer this
question for $\binom{\Omega}{\Tau}$; we only have a related result.

The following notion was introduced in \cite{tautau} as an
approximation for the notion of $\tau$-cover.
A family $Y\sbst\inf$ is \emph{linearly refinable}
if for each $y\in Y$ there exists an infinite subset
$\hat y\sbst y$ such that the family $\hat Y = \{\hat y : y\in Y\}$ is
linearly ordered by $\as$.
A cover $\cU$ of $X$ is a \emph{$\tau^*$-cover} of $X$ if
and $h_\cU[X]$ (where $h_\cU$ is the function defined before Lemma \ref{notions})
is linearly refinable.
By Lemma \ref{notions}, every $\tau^*$-cover is an $\omega$-cover, and
any $\tau$-cover is a $\tau^*$-cover.
Let $\Tau^*$, $\BTstar$, and $\CTstar$ denote the collections of all
\emph{countable} open, Borel, and clopen $\tau^*$-covers,
respectively.
\begin{thm}\label{OmChTau}
The property $\binom{\Omega}{\Tau^*}$ is
closed under taking finite powers.
\end{thm}
\begin{proof}
Fix $k$.
In \cite{coc2} it is proved that for each open $\omega$-cover $\cU$ of
$X^k$ there exists an open $\omega$-cover $\cV$ of $X$
such that the $\omega$-cover $\cV^k=\{V^k : V\in\cV\}$ of $X^k$ refines $\cU$.

Assume that $\cU$ is an open $\omega$-cover of
$X^k$. Choose an open $\omega$-cover $\cV$ of $X$
such that $\cV^k$ refines $\cU$.
Apply $\binom{\Omega}{\Tau^*}$ to choose a subcover
$\cW$ of $\cV$ such that $\cW$ is a $\tau^*$-cover of $X$.
Then $\cW^k$ is a $\tau^*$-cover of $X^k$ \cite{tautau}.
For each $W\in\cW$ choose $U_W\in\cU$ such that $W^k\sbst U_W$.
As $\tau^*$-covers are closed under taking de-refinements \cite{tautau},
$\{U_W : W\in\cW\}$ is a $\tau^*$-cover of $X$.
\end{proof}

Thus, if $\split(\Tau^*,\Tau^*)$ is closed under taking finite powers,
then so is
$\binom{\Omega}{\Tau^*}\cap\split(\Tau^*,\Tau^*) = \split(\Omega,\Tau^*)$.

We can get very close to showing that no class between
$\split(\BL,\BL)$ and $\split(\CT,\CT)$ is closed
under taking finite powers.

\begin{thm}
Assume \CH{} (or just $\t=\c$). Then there exist
sets of reals $L_0$ and $L_1$ such that:
\be
\i $L_0$ and $L_1$ satisfy $\split(\BO,\BO)$ and $\split(\BL,\BL)$,
\i $L = L_0\cup L_1$ satisfies $\split(\BL,\BL)$,
\i $L_0\x L_1$ and $L\x L$ do not satisfy $\split(\CTstar,\CL)$; and
\i $L_0\x L_1$ (and therefore $L\x L$) is not hereditarily-$\split(\CT,\CT)$.
\ee
In particular, the classes $\split(\Lambda,\Lambda)$ and $\split(\Omega,\Lambda)$
(and their Borel and clopen versions) are not closed under taking finite powers,
and $\split(\Omega,\Omega)$ (and its Borel and clopen versions) is not closed under taking finite products.
\end{thm}
\begin{proof}
The essence of the proof is the following lemma.
\begin{lem}\label{Lsets}
Assume \CH{} (or just $\t=\c$).
Then there exist $\t$-Luzin subsets $L_0=\{a^0_\alpha : \alpha<\c\}$ and $L_1=\{a^1_\alpha : \alpha<\c\}$
of $\inf$ such that $L_0$ and $L_1$ satisfy $\sone(\BO,\BO)$, and
$B = \{a^0_\alpha\cap a^1_\alpha : \alpha<\c\}$
is a simple $P$-point base.
\end{lem}
\begin{proof}
As we assume that $\t=\c$, there exists a simple $P$-point $U=\{a_\alpha : \alpha<\c\}$
(see the discussion before Corollary \ref{ShBlmodel}).

As $\t\le\add(\M)$, we have that $\add(\M)=\c$ and we can
repeat the construction given in \ref{notadd}, with the following modification:
At step $\alpha$ of the construction, consider the subset
$Y=\{a^0_\beta\cap a^1_\beta : \beta<\alpha\}\cup\{a_\alpha\}$ of $U$.
As $\alpha<\u$, this is not a base for $U$ and as $U$ is a simple
$P$-point, there exists $\tilde a_\alpha\in U$ which is a \pi{} of $Y$.
Now find, as done there, elements
$a^0_\alpha,a^1_\alpha\in U\sm Y_\alpha^*$ such that
$a^0_\alpha,a^1_\alpha\in (U\cap G_k)\sm Y_\alpha$,
and $a^0_\alpha \cap a^1_\alpha \as \tilde a_\alpha$.
\end{proof}

\begin{lem}\label{contcap}
The mapping from $\inf\x \inf$ to $\inf$ defined by
$$(a,b)\mapsto a\cap b$$
is continuous.
\end{lem}
\begin{proof}
It is enough to show that the preimage of a subbasic open set is
open.
Indeed, for each $n$ the preimage of $O_n = \{a\in\inf : n\in a\}$ is $O_n\x O_n$,
and the preimage of $\N\sm O_n$ is the union of the open sets
$O_n\x (\N\sm O_n)$, $(\N\sm O_n)\x O_n$, and $(\N\sm O_n)\x (\N\sm O_n)$.
\end{proof}

Let $U$, $L_0$, and $L_1$ be as in Lemma \ref{Lsets}.
By Lemma \ref{soneBOBO}, (1) holds.
As $L=L_0\cup L_1$ is an $\add(\M)$-Luzin set, (2) holds.
By Lemma \ref{contcap},
$B=\{a^0_\alpha\cap a^1_\alpha : \alpha<\c\}$ is a continuous image
of the subset $\Delta = \{(a^0_\alpha,a^1_\alpha) : \alpha<\c\}$ of $L_0\x L_1$.
As $B$ is a simple $P$-point base, we have by Lemma \ref{images2}
that $\Delta$ does not satisfy $\split(\CT,\CT)$. This proves (4).

To prove (3), we need to extend Lemma \ref{morenotions}.
Note that a base for a simple $P$-point need not be linearly ordered by $\as$, and
therefore need not be a simple $P$-point base according to our usage of this term.
\begin{lem}
Assume that $\cU=\{U_n\}_{n\in\N}$ is a cover of $X$.
The following are equivalent:
\be
\i $\cU$ is a $\tau^*$-cover of $X$ which cannot be split into two large covers of $X$; and
\i $h_\cU[X]$ is a base for a simple $P$-point.
\ee
\end{lem}
\begin{proof}
$1\Impl 2$: $\cU$ is, in particular, an $\omega$-cover which cannot be split into two
large covers.
By Lemma \ref{morenotions}, $Y=h_\cU[X]$ is base for a nonprincipal ultrafilter $U$ on $\N$.
By the definition of $\tau^*$-covers, $Y$ is linearly refinable.
Let $\hat Y$ be a linear refinement of $Y$. Then also $\hat Y$ is reaping, and clearly
it is centered. Thus, $\hat Y$ generates a nonprincipal filter $\tilde U$ containing
$U$. As $U$ is maximal, $U=\tilde U$ and $\hat Y$ witnesses that $U$ is a simple $P$-point.

$2\Impl 1$: Assume that $Y=h_\cU[X]$ is a base for a simple $P$-point $U$.
Choose a linearly ordered base $\hat Y$ for $U$.
Then for each $y\in Y$ there exists $\hat y\in\hat Y$ such that $\hat y\as y$.
Thus $\hat Y$ witnesses that $Y$ is linearly refinable.
\end{proof}
Consequently, a set of reals $X$ satisfies $\split(\CTstar,\CL)$ if, and only if,
every continuous image of $X$ in $\inf$ is not a base for a simple $P$-point.%
\footnote{Here too, the analogue Borel version also holds.
Moreover, we can show in a similar manner that
the combinatorial counterpart of $\lnot\split(\CTstar,\CO)$
and its Borel version is a \emph{subbase} for a simple $P$-point.}
This proves (3).
\end{proof}

With regard to finite products, only two problems remain open.
It seems that we will not take a great risk by stating them as a conjecture.
\begin{conj}\label{conj}
None of the classes $\split(\Tau,\Tau)$ and
$\binom{\Tau}{\Gamma}$ is provably closed under taking finite products.
\end{conj}

In the case of finite powers, we have more problems waiting for a solution.
\begin{prob}\label{powcl}
Is any of the classes $\split(\Omega,\Omega)$, $\binom{\Omega}{\Tau}\cap\split(\Tau,\Tau)$, or
$\split(\Tau,\Tau)$ closed under taking finite powers?
\end{prob}
The best candidate (if any) for a positive answer seems to be
$\split(\Omega,\Omega)$. Observe that the methods of \cite{coc2}
only give that if $X$ satisfies $\split(\Omega,\Omega)$, then
for each open $\omega$-cover $\cU$ of $X^k$
there exists a \emph{refinement} $\cV$ of $\cU$ such that
$\cV$ is an open $\omega$-cover of $X^k$ that can be
split into two disjoint $\omega$-covers of $X^k$.

We conclude this paper with the following related result.
As we mentioned in the introduction, it is proved in \cite{coc7}
that if all finite powers of $X$ have the Hurewicz property,
then $X$ satisfies $\split(\Omega,\Omega)$.
As the critical cardinality of the Hurewicz property is $\b$ and it is consistent that
$\b<\r$, the Hurewicz property is strictly stronger than $\split(\Lambda,\Lambda)$ \cite{coc2}.
Thus, the following theorem is strictly stronger than the quoted result.
\begin{thm}
Assume that for each $k$ $X^k$ satisfies $\split(\Omega,\Lambda)$.
Then $X$ satisfies $\split(\Omega,\Omega)$.
(The analogue assertions for the clopen and Borel cases also hold.)
\end{thm}
\begin{proof}
We say that $\cU$ is a \emph{$k$-cover} of $X$ if
($X$ is not contained in any member of $\cU$, and)
each $k$-element subset of $X$ is covered by some member of $\cU$.
Thus $\cU$ is a $k$-cover of $X$ if, and only if,
$$\cU^k = \{U^k : U\in\cU\}$$
is a cover of $X^k$.
Also, observe that $\cU$ is an $\omega$-cover of $X$ if, and only if,
$\cU^k$ is an $\omega$-cover of $X^k$.

\begin{lem}\label{Om-k}
Assume that $X^k$ satisfies $\split(\Omega,\Lambda)$.
Then each open $\omega$-cover $\cU$ of $X$ can be split
into two disjoint subsets $\cV$ and $\cW$
such that $\cV$ is an $\omega$-cover of $X$ and
$\cW$ is a $k$-cover of $X$.
\end{lem}
\begin{proof}
Assume that $\cU$ is an open $\omega$-cover of $X$.
Then for each $k$, $\cU^k$ is an $\omega$-cover of $X^k$,
and, by the assumption, can be split into two disjoint large covers
$\cV^k$ and $\cW^k$. Consequently, $\cV$ and $\cW$
are (large) $k$-covers of $X$.
As $\cU=\cV\cup\cW$ and the property of being an $\omega$-cover
is Ramseyan, at least one of the pieces
$\cV$ or $\cW$ is an $\omega$-cover of $X$.
\end{proof}
Assume that $\cU$ is an open $\omega$-cover of $X$.
As $X^2$ satisfies $\split(\Omega,\Lambda)$, we have
by Lemma \ref{Om-k} that $\cU=\cV_1\uplus\cW_1$
($\uplus$ denotes disjoint union) where
$\cV_1$ is an $\omega$-cover of $X$ and $\cW_1$
is a $2$-cover of $X$.
Continue inductively: Given an open $\omega$-cover $\cV_{k-1}$ ($k>1$)
of $X$, use the fact that $X^{k-1}$ satisfies $\split(\Omega,\Lambda)$ and
Lemma \ref{Om-k} to split $\cV_{k-1} = \cV_k\uplus\cW_k$ such
that $\cV_k$ is an $\omega$-cover of $X$ and $\cW_k$ is an $k+1$-cover of $X$.
Set
$$\cU_1 = \Union_{n\in\N}\cW_{2n+1},\ \cU_2 = \Union_{n\in\N}\cW_{2n}.$$
Then $\cU_1$ and $\cU_2$ are disjoint subcovers of $\cU$, and they are
$k$-covers of $X$ for all $k$, that is, $\omega$-covers of $X$.
\end{proof}
Thus, in order to prove that $\split(\Omega,\Omega)$ is closed under
taking finite powers, it is enough to show that all finite powers of
members of $\split(\Omega,\Omega)$ satisfy $\split(\Omega,\Lambda)$.

\section{Summary of open problems}
One may argue that the property $\split(\fU,\fV)$ is only (or, at least, more)
interesting when $\fU\sbst\fV$. If we accept this thesis, then no classification
problem (Part 1) remains open, and the more interesting problems in Part 2 are
Problems \ref{SpLamLamAdd}, \ref{BorelHered} (for the first three properties),
\ref{conj} (for the first property),
and \ref{powcl} (for the first and last properties).

On the other hand, the other problems
(\ref{survprobs}, \ref{BorelHered} for the fourth property,
\ref{conj} for the second property, and \ref{powcl}
for the second property),
which involve properties of the form
$\binom{\fU}{\fV}$, rise naturally in many other contexts, published (e.g, \cite{tau, tautau, huremen2, ideals})
and unpublished. In this sense, these problems are not less, and maybe more, interesting.

\bigskip
\paragraph{\textbf{Acknowledgements.}}
We thank Andreas Blass for the useful details on reference \cite{dordal},
Saharon Shelah for the proof of Lemma \ref{UminusM}, and the referee for
the extension of Theorem \ref{add3} to its current form.

\end{document}